\documentclass[12pt]{amsart}
\usepackage[mathscr]{eucal}
\usepackage{amsthm}
\usepackage{amssymb}
\usepackage{amscd}
\usepackage{verbatim}

\usepackage{array}
\usepackage{supertabular}
\usepackage{rotating}

\makeatletter
\@ifclasslater{amsart}{1999/11/24}{}{
\renewcommand*\subjclass[2][1991]{%
  \def\@subjclass{#2}%
  \@ifundefined{subjclassname@#1}{%
    \ClassWarning{\@classname}{Unknown edition (#1) of Mathematics
      Subject Classification; using '1991'.}%
  }{%
    \@xp\let\@xp\subjclassname\csname subjclassname@#1\endcsname
  }%
}
\renewcommand{\subjclassname}{%
  \textup{1991} Mathematics Subject Classification}
\@xp\let\csname subjclassname@1991\endcsname \subjclassname
\@namedef{subjclassname@2000}{%
  \textup{2000} Mathematics Subject Classification}
}
\makeatother

\def\Bbb{\mathbb}
\def\frak{\mathfrak}

\newenvironment{pf*}[1]{\proof[#1]}{\endproof}
\newtheorem{Fact}[equation]{Fact}

\theoremstyle{definition}
\newtheorem{Definition}[equation]{Definition}

\theoremstyle{remark}
\newtheorem{Remark}[equation]{Remark}

\numberwithin{equation}{section}

\newcommand{\factref}[1]{Fact~\ref{#1}}

\newcommand{\defeq}{\overset{\operatorname{\scriptstyle def.}}{=}}
\newcommand{\C}{{\Bbb C}}
\newcommand{\Z}{{\Bbb Z}}

\newcommand{\g}{{\frak g}}

\newcommand{\ve}{\varepsilon}
\newcommand{\Uq}{{\mathbf U}_q(\mathfrak g)} 

\newcommand{\Ua}{{\mathbf U}_q(\widehat{\mathfrak g})} 
\newcommand{\Ul}{{\mathbf U}_q({\mathbf L}{\mathfrak g})} 

\newcommand{\bfR}{\mathbf R}
\newcommand{\Lg}{\mathbf L\g}
\newcommand{\qch}[1]{\widetilde{\chi_{#1}}}

\newcommand{\vi}{\varpi_5}
\newcommand{\vii}{\varpi_3+\varpi_7}
\newcommand{\viii}{\varpi_2+\varpi_8}
\newcommand{\viv}{\varpi_1+2\varpi_7}
\newcommand{\vv}{\varpi_1+\varpi_6}
\newcommand{\vvi}{2\varpi_2}
\newcommand{\vvii}{\varpi_7+\varpi_8}
\newcommand{\vviii}{\varpi_1+\varpi_3}
\newcommand{\vix}{\varpi_4}
\newcommand{\vx}{2\varpi_1+\varpi_7}
\newcommand{\vxi}{\varpi_2+\varpi_7}
\newcommand{\vxii}{\varpi_1+\varpi_8}
\newcommand{\vxiii}{2\varpi_7}
\newcommand{\vxiv}{\varpi_6}
\newcommand{\vxv}{3\varpi_1}
\newcommand{\vxvi}{\varpi_1+\varpi_2}
\newcommand{\vxvii}{\varpi_3}
\newcommand{\vxviii}{\varpi_1+\varpi_7}
\newcommand{\vxix}{\varpi_8}
\newcommand{\vxx}{2\varpi_1}
\newcommand{\vxxi}{\varpi_2}
\newcommand{\vxxii}{\varpi_7}
\newcommand{\vxxiii}{\varpi_1}
\newcommand{\vxxiv}{0}

\begin{document}
\title[$t$--analogs of $q$--characters of type $E_n$]
{$t$--analogs of $q$--characters
\\
of quantum affine algebras of type $E_6$, $E_7$, $E_8$
}
\author{Hiraku Nakajima}
\address{Department of Mathematics, Kyoto University, Kyoto 606-8502,
Japan
}
\email{nakajima@kusm.kyoto-u.ac.jp}
\urladdr{http://www.math.kyoto-u.ac.jp/\textasciitilde nakajima}
\thanks{Supported by the Grant-in-aid
for Scientific Research (No.17340005), the Ministry of Education,
Japan.}
%
\subjclass[2000]{Primary 17B37;
Secondary 14D21, 14L30, 16G20}
\begin{abstract}
We compute $t$--analogs of $q$--characters of all $l$--fundamental
representations of the quantum affine algebras of type $E_6^{(1)}$,
$E_7^{(1)}$, $E_8^{(1)}$ by a supercomputer. In particular, we prove the
fermionic formula for Kirillov-Reshetikhin modules conjectured by
Hatayama et al.\ \cite{HKOTY} for these classes of representations. We
also give explicitly the monomial realization of the crystal of the
corresponding fundamental representations of the qunatum enveloping
algebras associated with finite dimensional Lie algebras of types
$E_6$, $E_7$, $E_8$. These are computations of Betti numbers of graded
quiver varieties, quiver varieties and determination of all
irreducible components of the lagrangian subvarities of quiver
varieties of types $E_6$, $E_7$, $E_8$ respectively.
\end{abstract}
\maketitle

\section*{Introduction}

Let $\g$ be a simple Lie algebra of type $ADE$ over $\C$ with the
index set $I$ of simple roots, $\Lg = \g \otimes \C[z,z^{-1}]$ be its
loop algebra, and $\Ul$ be its quantum universal enveloping algebra,
or the quantum loop algebra for short.  It is a subquotient of the
quantum affine algebra $\Ua$, i.e.\ without central extension and
degree operator. It contains the quantum enveloping algebra $\Uq$
associated with $\g$ as a subalgebra.

By Drinfeld \cite{Dr} and Chari-Pressley \cite{CP-rep}, simple
$\Ul$-modules are parametrized by $I$-tuples of polynomials $P =
(P_i(u))_{i\in I}$ with normalization $P_i(0) = 1$. They are called
{\it Drinfeld polynomials}. Let us denote by $L(P)$ the simple module
with Drinfeld polynomial $P$. When $P$ is given by $P_i(u) =
(1-au)^{\delta_{iN}}$ for a given $N\in I$, we call corresponding
module an {\it $N^{\mathrm{th}}$ $l$--fundamental representation}.
(It has been called a {\it level $0$ fundamental module\/} or simply
{\it fundamental representation\/} in some literature.) We can assume
$a=1$ without the loss of generality as the general module is a
pullback of the module with $a=1$ by an algebra automorphism of $\Ul$.

Let $\chi_{q,t}(L(P))$ be the $t$--analog of $q$--character of a
simple module $L(P)$ defined by the author
\cite{Na-qchar,Na-qchar-main}. It is defined via the geometry of
graded quiver varieties.
It values in certain Laurent polynomial ring with infinitely many
variables with integer coefficients. It is a $t$--analog of the
$q$--character $\chi_q(L(P))$ introduced earlier \cite{Kn,FR}, which
was a refinement of the ordinary character of the restriction of
$L(P)$ to a $\Uq$-module.
In \cite{Na-qchar,Na-qchar-main} we ``computed'' $\chi_{q,t}(L(P))$
for arbitrary given $L(P)$, in the sense that we gave a purely combinatorial
algorithm to write down all monomials and coefficients in
$\chi_{q,t}(L(P))$ where the final expression involves only $+$, $\times$,
integers and variables.

In order to clarify in what sense our result is new compared with
earlier results, we {\it define\/} what the word compute mean
precisely.
When we write the word compute in the quotation marks, it means that
we give a combinatorial algorithm to compute something in the above sense.
It does not necessarily mean that we actually compute it. We can
write a computer program in principle, but the question whether we can
actually compute it or not depends on the size of computer memory.
(For example, it is clear that the rank $n$ of $\g$ cannot be larger
than the size of the memory.)
On the other hand, when we write the word compute without the
quotation mark, we mean to compute something in a strict sense, i.e.\
we express something so that it contains only finitely many $\pm$,
$\times$, integers and variables.
For example, if we write $x = \sum_{i=1}^{2^{(2^{100})}} a_i$ for some explicit
$a_i$, we ``compute'' $x$, but we do not compute $x$ unless we actually
compute the sum.
On the other hand, we do not require that the final expression can be
read by the human, as such a concept cannot make precise.

The algorithm is separated into three steps:
\begin{enumerate}
\item ``Computation'' of $\chi_{q,t}$ for $l$--fundamental representations.
\item ``Computation'' of $\chi_{q,t}$ for standard modules, i.e.\ tensor
  products of $l$--fundamental representations.
\item ``Computation'' of the $t$-analog of the composition factors of
simple modules in standard modules.
\end{enumerate}
The third step is analogous to the definition of Kazhdan-Lusztig
basis. If $M(P)$ denote the standard module, we have
\begin{equation}\label{eq:L(P)}
    \overline{\chi_{q,t}(L(P))} = \chi_{q,t}(L(P)), \qquad
    \chi_{q,t}(L(P)) = \chi_{q,t}(M(P)) + \sum_{Q:Q<P} a_{PQ}(t) 
    \chi_{q,t}(M(Q))
\end{equation}
for some $a_{PQ}(t)\in t^{-1}\Z[t^{-1}]$, where `$<$' is a certain 
explicitly defined ordering. Thus $a_{PQ}(t)$ is analogous to
Kazhdan-Lusztig polynomials. The above characterization allows us to
``compute'' $a_{PQ}(t)$, once $\chi_{q,t}(M(P))$ is ``computed''. (And
it is known that the actual computation of Kazhdan-Lusztig polynomials
is very hard.)

In the second step, we express $\chi_{q,t}(M(P))$ as a twisted
multiplication of $\chi_{q,t}$ of $l$--fundamental representations. It
is almost the same as usual multiplication on the polynomials, but a
product of two monomials $m$, $m'$ is twisted as $t^{2d(m,m')} mm'$.
Therefore this step is very simple. It is clear that
$\chi_{q,t}(M(P))$ can be ``computed'' if $\chi_{q,t}$ of
$l$--fundamental representations are ``computed''.

This paper concerns the first step. Our ``computation'' in
\cite{Na-qchar,Na-qchar-main} was $t$--analog of the ``computation'' by
Frenkel-Mukhin~\cite{FM}. 
It is based on the observation that (a) $\chi_{q,t}$ satisfies a
certain analog of the Weyl group invariance of the ordinary
characters, and (b) the $l$--fundamental representation satisfies a
certain property analogous to that of minuscule representations of
$\g$. Recall that a simple finite dimensional representation of $\g$
is called {\it minuscule\/} if all weights are conjugates of the
highest weight under the Weyl group each occurring with multiplicity
$1$.

When $\g$ is of classical type, i.e.\ of type $A$, $D$, the author
gave a tableaux sum expression of $\chi_{q,t}$ of $l$--fundamental
representations \cite{Na-AD}. It means that we give another
``computation'' of $\chi_{q,t}$, which are more familiar to us than
the above one. It does not mean we compute $\chi_{q,t}$ in our
strict sense. In fact, the comparison of two methods does not make sense
unless we define what we mean by `familiar'.
In practice, it just means that we have a faster algorithm for the
actual computer calculation.

In this paper we report the actual computer computation of
$\chi_{q,t}$ of $l$--fundamental representations when $\g$ is of type
$E_6^{(1)}$, $E_7^{(1)}$, $E_8^{(1)}$. Our algorithm is implemented in
the computer language {\bf C}. The source code is available at
\linebreak[3] \verb|http://www.math.kyoto-u.ac.jp/~nakajima/Qchar/|.
The author's personal computer (Dell Dimension 9100) can give the
answer up to the $6^{\mathrm{th}}$ $l$--fundamental representation of
$E_8$, where our numbering of $I$ is the following:
\begin{equation*}
\begin{array}{ccccccccccccc}
7 & - & 6 & - & 5 & - & 4 & - & 3 & - & 2 & - & 1
\\
&&&& | &&&&&&&&
\\
&&&& 8 &&&&&&&&
\end{array}  
\end{equation*}
We need about 120Mbtyes of the memory for this calculation. For the
$4^{\mathrm{th}}$ and $5^{\mathrm{th}}$ $l$--fundamental
representations, the computation was done on a supercomputer FUJITSU
HPC 2500 at Kyoto University. The calculation required about 2.6Gbytes
(for $4^{\mathrm{th}}$) and 120Gbytes (for $5^{\mathrm{th}}$) of
memory, and it took 6 hours and 350 hours for the calculation
respectively. The final answers (stored in a compressed format as
explained below) are 3.2Gbytes and 180Gbytes respectively.
In fact, the calculation of the $4^{\mathrm{th}}$ one was done several
years ago and was mentioned in some of the author's papers. However
we needed to wait for the Kyoto University to renovate the
supercomputer so that we can use 120Gbytes of memory in a single
program, and then wait for the author to get an enough budget to use
the supercomputer.

As far as the author knows, the computation (in our strict sense) for
the $5^{\mathrm{th}}$ one was not known before.
Frenkel-Mukhin, Hernandez-Schedler told the author that they wrote
computer programs calculating $\chi_{q,t=1}$ and $\chi_{q,t}$ respectively.
But both had a problem of computer memory. 

In conclusion, we can now delete the quotation mark for computation in
the first step of the algorithm for type $E$ above.

As an application, we can compute $t$--analog of the ordinary
characters of the restrictions of $l$--fundamental representations to
$\Uq$-modules. The $l$--fundamental modules are examples of the
so-called Kirillov-Reshetikhin modules. Kirillov-Reshetikhin gave
conjectural formula for the ordinary character of the restriction of a
Kirillov-Reshetikhin module \cite{KR}. Its graded version (i.e.\
$t$--analog) together with an interpretation in terms of the
conjectural crystal base was given by Hatayama, Kuniba, Okado, Takagi
and Yamada \cite{HKOTY}. Then Lusztig conjectured that their
conjectural grading is the same as the cohomological degree
\cite{Lu:ferm}, in a certain class of Kirillov-Reshetikhin modules
including $l$--fundamental representations. Therefore the formula in
\cite{HKOTY}, in the class, gives the generating function of
Poincar\'e polynomials of quiver varieties.
In general, the conjectural formula is expressed as a summation over
partitions, and called a {\it fermionic formula}.
The author gave an expression for $t=1$ in \cite[Cor.~1.3]{Na-KR} (the
result was extended to type $BCFG$ in \cite{Her}). It is again given
as a summation over partition, but the definition of the binomial
coefficient appearing in the coefficients is different. The
equivalence between two expressions are not known so far, therefore
the original fermionic formula is remained open.

For an $l$--fundamental representation, the original fermionic formula
can be given by an explicit polynomial by the so-called Kleber's
algorithm \cite{Kl}. Here we do not make precise what we mean by
`explicit'.
For types $A$, $D$, it was shown in \cite{Na-AD} that this
`explicit' expression for an $l$--fundamental representation is equal
to the ``computation'' in \cite{Na-qchar}. For type $E$, the algorithm
can be used to compute the fermionic formula in our strict sense. Then
the result can be checked in some special cases previously computed
(at least for $t=1$) (e.g.\ \cite{CP}), but most of $l$--fundamental
representations have remained open.
Remark that Kleber's algorithm does not apply to the modified formula
in \cite{Na-KR}, so it is not known that the modified formula gives
the computation in the strict sense.

Our computation of $\chi_{q,t}$ gives the explicit expression and we
find that it is the same as one given in \cite{HKOTY}.
Therefore we prove Lusztig's conjecture for all
$l$--fundamental representations.

Also as another application, we determine all monomials appearing
in the monomoial realization of the crystal corresponding to
fundamental representations of type $E$. For types $A$, $D$, they were
determined in \cite{Na-AD} as an application of the explicit
description of $\chi_{q,t}$ of $l$--fundamental representations. 
For types $B$, $C$, they were determined in \cite{kks}. For types $F$,
$G$, they can be easily determined (cf.\ \cite{HN}). In conclusion, we
describe the monomial realization of the crystals of all
fundamental representations explicitly.

\subsection*{Acknowledgement}

A part of the computer program was written while the author stayed at
Centre for Advanced Study (CAS) at the Norwegian Academy of Science
and Letters in 2002. He would like to thank CAS for the hospitality.

\section{$t$--analogs of $q$--characters}

We shall not give the definition of quantum loop algebras, nor
their finite dimensional representations in this paper. (See
\cite{Na-qchar} for a survey.) We just review properties of
$\chi_{q,t}$, as axiomized in \cite{Na-qchar-main}.


Let
\(
\mathscr Y_t \defeq 
   \Z[t,t^{-1},Y_{i,a}, Y_{i,a}^{-1}]_{i\in I, a\in\C^*}
\)
be a Laurent polynomial ring of uncontably many variables $Y_{i,a}$'s
with coefficients in $\Z[t,t^{-1}]$. A {\it monomial\/} in $\mathscr
Y_t$ means a monomial only in $Y_{i,a}^\pm$, containing no
$t$'s. Therefore a polynomial is a sum of monomials multiplied by
Laurent polynomials in $t$, called coefficients as usual.
Let
\begin{equation*}
   A_{i,a} \defeq Y_{i,a\ve} Y_{i,a\ve^{-1}}
     \prod_{j:j\neq i} Y_{j,a}^{c_{ij}},
\end{equation*}
where $c_{ij}$ is the $(i,j)$-entry of the Cartan matrix.
Let $\mathcal M$ be the set of monomials in $\mathscr Y_t$.

\begin{Definition}
(1) For a monomial $m\in\mathcal M$, we define $u_{i,a}(m)\in\Z$ be the
degree in $Y_{i,a}$, i.e.\ 
\begin{equation*}
   m = \prod_{i,a} Y_{i,a}^{u_{i,a}(m)}.
\end{equation*}

(2) A monomial $m\in\mathcal M$ is said {\it $i$--dominant\/} if
$u_{i,a}(m)\ge 0$ for all $a$. It is said {\it l--dominant\/} if it is
$i$--dominant for all $i$.

(3) Let $m, m'$ be monomials in $\mathcal M$. We say $m \le m'$ if
$m/m'$ is a monomial in $A_{i,a}^{-1}$ ($i\in I$, $a\in\C^*$).
Here a monomial in $A_{i,a}^{-1}$ means a product of nonnegative
powers of $A_{i,a}^{-1}$. It does not contain any factors
$A_{i,a}$. In such a case we define $v_{i,a}(m, m')\in\Z_{\ge 0}$ by
\begin{equation*}
   m = m' \prod_{i,a} A_{i,a}^{-v_{i,a}(m,m')}.
\end{equation*}
This is well-defined since the $\ve$-analog of the Cartan matrix is
invertible. We say $m < m'$ if $m\le m'$ and $m\neq m'$.

(4) For an $i$--dominant monomial $m\in\mathcal M$ we define
\begin{equation*}
   E_i(m) \defeq
    m\, \prod_a
     \sum_{r_a=0}^{u_{i,a}(m)}
     t^{r_a(u_{i,a}(m)-r_a)}
     \begin{bmatrix}
       u_{i,a}(m) \\ r_a
     \end{bmatrix}_t A_{i,a\ve}^{-r_a},
\end{equation*}
where
\(
\left[\begin{smallmatrix}
  n \\ r
\end{smallmatrix}\right]_t
\)
is the $t$-binomial coefficient.

(5) We define a ring involution
$\setbox5=\hbox{A}\overline{\rule{0mm}{\ht5}\hspace*{\wd5}}$ on
${\mathscr Y}_t$ by $\overline{t} = t^{-1}$,
$\overline{Y_{i,a}^\pm} = Y_{i,a}^\pm$.
\end{Definition}

Suppose that {\it l\/}--dominant monomials $m_{P^1}$, $m_{P^2}$ and
monomials $m^1\le m_{P^1}$, $m^2\le m_{P^2}$ are given. We define an
integer $d(m^1, m_{P^1}; m^2, m_{P^2})$ by
\begin{multline}\label{eq:d}
   d(m^1, m_{P^1}; m^2, m_{P^2})
\\
   \defeq
   \sum_{i,a} \left( v_{i,a\ve}(m^1, m_{P^1}) u_{i,a}(m^2)
   + u_{i,a\ve}(m_{P^1}) v_{i,a}(m^2, m_{P^2})\right).
\end{multline}

For an $I$-tuple of rational functions $Q/R = (Q_i(u)/R_i(u))_{i\in
I}$ with $Q_i(0) = R_i(0) = 1$, we set
\begin{equation*}
   m_{Q/R} \defeq
   \prod_{i\in I} \prod_{\alpha} \prod_{\beta}
     Y_{i,\alpha} Y_{i,\beta}^{-1},  
\end{equation*}
where $\alpha$ (resp.\ $\beta$) runs roots of $Q_i(1/u) = 0$
(resp.\ $R_i(1/u) = 0$), i.e.\ 
$Q_i(u) = \prod_\alpha ( 1 - \alpha u)$ (resp.\ $R_i(u) = \prod_\beta
(1 - \beta u)$). As a special case, an $I$-tuple of polynomials $P =
(P_i(u))_{i\in I}$ defines $m_P = m_{P/1}$. The $l$--dominant monomial
$m_{P^\alpha}$ appeared above is assoicated to an $I$-tuple of
polynomials $P = (P_i(u))_{i\in I}$.
In this way, the set $\mathcal M$ of monomials are identified with the 
set of $I$-tuple of rational functions, and the set of {\it
l\/}--dominant monomials are identified with the set of $I$-tuple of
polynomials. 

The $t$--analog of the Grothendieck ring $\bfR_t$ 
is a free $\Z[t,t^{-1}]$-module with base $\{ M(P) \}$ where $P =
(P_i(u))_{i\in I}$ is the Drinfeld polynomial.
(We do not recall the definition of standard modules $M(P)$ here, but
the reader safely consider them as formal variables.)

The $t$--analog of the $\ve$--character homomorphism is a
$\Z[t,t^{-1}]$-linear homomorphism
\(
   \chi_{q,t}\colon \bfR_t \to \mathscr Y_t.
\)
It is defined as the generating function of Poincar\'e polynomials of
graded quiver varieties, or the generating function of graded dimensions of
$l$--weight spaces of a $\Ul$-module \cite{VV2}, and will not be
reviewed in this paper. 

We also need a slightly modified version:
\begin{equation*}
   \qch{q,t}(M(P)) = \sum_m t^{d(m,m_P;m,m_P)} a_m(t) m
   \qquad \text{if $\chi_{q,t}(M(P)) = \sum_m a_m(t) m$}.
\end{equation*}
If we know one of $\chi_{q,t}$ and $\qch{q,t}$, we know the remaining
one.

The following was proved in \cite{Na-qchar,Na-qchar-main}:
\begin{Fact}\label{thm:ind}
\textup{(1)} 
The $\chi_{q,t}$ of a standard module $M(P)$ has a form
\begin{equation*}
   \chi_{q,t}(M(P)) = m_P + \sum a_m(t) m,
\end{equation*}
where the summation runs over monomials $m < m_P$.

\textup{(2)}
For each $i\in I$, $\qch{q,t}(M(P))$ can be expressed as a linear
combination \textup(over $\Z[t,t^{-1}]$\textup) of $E_i(m)$ with
$i$--dominant monomials $m$.

\textup{(3)}
Suppose that two $I$-tuples of polynomials $P^1 = (P^1_i)$, $P^2 =
(P^2_i)$ satisfy the following condition:
\begin{equation}
\label{eq:Z}
\begin{minipage}[m]{0.75\textwidth}
\noindent   
   $a/b\notin\{ \ve^n \mid n\in\Z, n \ge 2\}$ for any
   pair $a$, $b$ with $P^1_i(1/a) = 0$, $P^2_j(1/b) =
   0$ \textup($i,j\in I$\textup).
\end{minipage}
\end{equation}
Then we have
\begin{equation*}
  \qch{q,t}(M(P^1P^2)) =
  \sum_{m^1, m^2} t^{2d(m^1, m_{P^1}; m^2, m_{P^2})}
   a_{m^1}(t) a_{m^2}(t) m^1 m^2
,
\end{equation*}
where
\(
   \qch{q,t}(M(P^a)) = \sum_{m^a} a_{m^a}(t) m^a
\)
with $a=1,2$.

Moreover, properties \textup{(1),(2),(3)} uniquely determine
$\chi_{q,t}(M(P))$.

\textup{(4)} The $\chi_{q,t}$ of the simple module $L(P)$ is given by 
\eqref{eq:L(P)}.
\end{Fact}

Apart from the existence problem, one can consider the above
properties (1), (2), (3) as the definition of $\chi_{q,t}$ (an
axiomatic definition). We only use the above properties, and the
reader can safely forget the original definition. Note that we will
prove the existence of $\chi_{q,t}$ by our computer calculation.

By the property (1) we call the monomial $m_P$ corresponding to the
Drinfeld polynomial $P$ {\it $l$--highest weight monomial}.

\section{Algorithm}\label{sec:algorithm}

In this section we shall explain our algorithm to determine
$\qch{q,t}(L(P))$ recursively starting from the $l$--dominant weight
monomial $m_P$. It is a slight modification of one in \cite{FM}.
We shall also explain why we require large memory to compute
$\chi_{q,t}$ of the $5^{\mathrm{th}}$ $l$--fundamental representation
of $\Ul$ with $\g = E_8$. The problem does not exist for the other
$l$--fundamental representations. 

We take a Drinfeld
polynomial $P = (P_i(u))$ $P_i(u) = (1-u)^{\delta_{iN}}$ corresponding
to the $N^{\mathrm{th}}$ $l$--fundamental representation.

One of the key property of $\chi_{q,t}$ of an $l$--fundamental
representation is that all monomials appearing in $\chi_{q,t}$ are not
$l$--dominant except the $l$--highest one. This was proved in
\cite[Cor.~4.5]{FM} and \cite[4.13]{Na-qchar-main}.

For each monomial $m$ in $\qch{q,t}(L(P))$ we determine the
coefficient $a_m(t)\in\Z[t]$ and the $I$-tuple of polynomial
$(a_{m,i}(t))_{i\in I}\in \Z[t]$ (called {\it coloring\/})
recursively. Let us introduce several concepts.
We say $m$ is {\it admissible\/} if all
$a_{m,i}(t)$ are the same for any $i$ such that $m$ is not
$i$--dominant. We say {\it the algorithm fails at $m$\/} if $m$ is not
admissible. 
We say {\it the algorithm stops at $m$} if $m$ is
$l$--dominant. 

Now we explain the algorithm.
At the first stage we set $a_{m_P}(t) = 1$ and $a_{m_P,i}(t) = 0$ for
all $i\in I$ for the $l$--highest weight monomial $m_P$.
Next take a monomial $m$ such that $a_m(t)$ and $a_{m,i}(t)$ are
determined. If $m$ is not $i$--dominant for any $i$ (this will happen
if $m$ the $l$--lowest weight vector), we do nothing on $m$ and go to
the next monomial.
If $m$ is $i$--dominant, we compute $(a_m(t) - a_{m,i}(t))
E_i(m)$. We call this procedure the {\it $i$-expansion at $m$}.
We add a monomial $m'$ appearing there to the list.
And for a monomial $m'$ in the list, we set $a_{m',i}(t)$ be the sum
of the contribution to $m'$ in the $i$-expansion at $m$ for various $m
< m'$ which is $i$--dominant.
As there is only finitely many $m < m'$, $a_{m',i}(t)$ will be
eventually determined. 
After all $a_{m',i}(t)$ are determined in this way, we can ask
$m'$ is admissible or not.
If $m'$ is not admissible (i.e.\ the algorithm fails at $m'$), we stop.
If $m'$ is $l$--dominant (i.e.\ the algorithm stop at $m'$), we stop.
If $m'$ is admissible and not $l$--dominant, we set $a_{m'}(t) =
a_{m',i}(t)$ for some (and any by admissibility) $i$ such that $m'$ is
not $i$--dominant.
We continue this procedure until all $a_{m}(t)$ and
$a_{m,i}(t)$ are determined, and all $(a_m(t) - a_{m,i}(t)) E_i(m)$
are expanded, or we stop at some $m$.

Now we apply the algorithm starting from the $l$--highest weight
monomial $m_P$.
As $\qch{q,t}(L(P))$ satisfies the properties (1),(2) in
\factref{thm:ind}, the algorithm cannot fail. As $\qch{q,t}(L(P))$ does
not contain $l$--dominant monomials other than $l$--highest one, the
algorithm cannot stop.
Finally as $L(P)$ is a finite dimensional, $\qch{q,t}(L(P))$ contains
only finitely many monomials. Therefore we eventually determine
all $a_{m}(t)$ and $a_{m,i}(t)$.

\begin{Remark}
If we apply the same algorithm in case $\g$ is a Kac-Moody Lie algebra
(say an affine Lie algebra), the algorithm does not fail, does not
stop, but we always get a new monomial in the expansion. Therefore the
procedure never end.
\end{Remark}

Now we consider the $5^{\mathrm{th}}$ $l$--fundamental representation
of $\Ul$ with $\g = E_8$ and we will explain the reason why we need
various tricks to save the size of data.
Because of these tricks, we had not known how big the total size is
in advance, so we used the following guess:
We know that the dimension of the $4^{\mathrm{th}}$ fundamental
representation of $\g$ is $146325270$, while $5^{\mathrm{th}}$ one is
$6899079264$. Therefore we expect that the corresponding
$\chi_{q,t}$'s have a similar ratio. We first compute the
$4^{\mathrm{th}}$ $l$--fundamental representations and expect that the
total size of the $5^{\mathrm{th}}$ one is about $50$ times as much.
This turned out to be approximately correct as we can see from the
data in Introduction.

By \cite[Prop.~3.4]{Na-AD} the set of monomials appearing in the
$q$--character of an $l$--fundametal representation has a
$\Uq$-crystal structure, which is isomorphic to the corresponding
fundamental representation of $\Uq$. In particular, the number of the
monomials appearing in the $5^{\mathrm{th}}$ $l$--fundamental
representation is equal to the dimension of the $5^{\mathrm{th}}$
fundamental representation of $\g = E_8$, i.e.\ $6899079264 \approx
6.4 \times 2^{30} = 6.4\mathrm{Giga}$. 
For each monomial $m$, we must remember (a) the expression of the
monomial and (b) the coloring, i.e.\ an $I$-tuple of polynomials in
$t$.

Let us first consider how we can express the monomial.
It is known that $l$--lowest weight monomial, i.e.\ the unique
monomial with (ordinary) weight $-\varpi_5$, is $Y_{5,q^{30}}^{-1}$
(see e.g.\ \cite[6.8]{FM}). We have
{\allowdisplaybreaks
\begin{equation*}
  \begin{split}
  Y_{5,q^{30}}^{-1} = Y_{5,1}
  & \times A_{1,q^{5}} A_{1,q^{7}} A_{1,q^{9}} A_{1,q^{11}} A_{1,q^{13}} A_{1,q^{15}}^2 A_{1,q^{17}} A_{1,q^{19}}
  A_{1,q^{21}} A_{1,q^{23}} A_{1,q^{25}}
\\ &\times
  A_{2,q^{4}} A_{2,q^{6}}^2 A_{2,q^{8}}^2 A_{2,q^{10}}^2 A_{2,q^{12}}^2 A_{2,q^{14}}^3 A_{2,q^{16}}^3 A_{2,q^{18}}^2
  A_{2,q^{20}}^2 A_{2,q^{22}}^2 A_{2,q^{24}}^2 A_{2,q^{26}}
\\ &\times
  A_{3,q^{3}} A_{3,q^{5}}^2 A_{3,q^{7}}^3 A_{3,q^{9}}^3 A_{3,q^{11}}^3 A_{3,q^{13}}^4 A_{3,q^{15}}^4 A_{3,q^{17}}^4 A_{3,q^{19}}^3
  A_{3,q^{21}}^3 A_{3,q^{23}}^3 A_{3,q^{25}}^2 A_{3,q^{27}}
\\ &\times
 A_{4,q^{2}} A_{4,q^{4}}^2 A_{4,q^{6}}^3 A_{4,q^{8}}^4 A_{4,q^{10}}^4 A_{4,q^{12}}^5 A_{4,q^{14}}^5 A_{4,q^{16}}^5 A_{4,q^{18}}^5
 A_{4,q^{20}}^4 A_{4,q^{22}}^4 A_{4,q^{24}}^3 A_{4,q^{26}}^2 A_{4,q^{28}}
\\ &\times
 A_{5,q^{1}} A_{5,q^{3}}^2 A_{5,q^{5}}^3 A_{5,q^{7}}^4 A_{5,q^{9}}^5 A_{5,q^{11}}^6 A_{5,q^{13}}^6 A_{5,q^{15}}^6 A_{5,q^{17}}^6
 A_{5,q^{19}}^6 A_{5,q^{21}}^5 A_{5,q^{23}}^4 A_{5,q^{25}}^3 A_{5,q^{27}}^2 A_{5,q^{29}}
\\ &\times
 A_{6,q^{2}} A_{6,q^{4}}^2 A_{6,q^{6}}^2 A_{6,q^{8}}^3 A_{6,q^{10}}^4 A_{6,q^{12}}^4 A_{6,q^{14}}^4 A_{6,q^{16}}^4 A_{6,q^{18}}^4
 A_{6,q^{20}}^4 A_{6,q^{22}}^3 A_{6,q^{24}}^2 A_{6,q^{26}}^2 A_{6,q^{28}}
\\ &\times
 A_{7,q^{3}} A_{7,q^{5}} A_{7,q^{7}} A_{7,q^{9}}^2 A_{7,q^{11}}^2 A_{7,q^{13}}^2 A_{7,q^{15}}^2 A_{7,q^{17}}^2 A_{7,q^{19}}^2
 A_{7,q^{21}}^2 A_{7,q^{23}} A_{7,q^{25}} A_{7,q^{27}}
\\ &\times
 A_{8,q^{2}} A_{8,q^{4}} A_{8,q^{6}}^2 A_{8,q^{8}}^2 A_{8,q^{10}}^3 A_{8,q^{12}}^3 A_{8,q^{14}}^3 A_{8,q^{16}}^3 A_{8,q^{18}}^3
 A_{8,q^{20}}^3 A_{8,q^{22}}^2 A_{8,q^{24}}^2 A_{8,q^{26}} A_{8,q^{28}}.
  \end{split}
\end{equation*}
Any} other monomial is given equal to $Y_{5,1}$ multiplied by a part of
$A_{i,q^k}$'s appeared above. We record the monomial as a sequence of
$A_{i,q^k}^m$'s, where $i$ runs $1$ to $8$, $k$ runs from $1$ to $29$,
and $m$ runs from $1$ to $6$. We can store the triple $(i, k, m)$ in a
single \verb+short int+, i.e.\ $16$bit of memory.
The length of the sequence is at most $106$, which is the length for
$Y_{5,q^{30}}^{-1}$. A naive count gives
\(
   6899079264 \times 106 \times 16\mathrm{bit} > 1300\mathrm{Gbyte}.
\)
This is too large. Therefore we use the following trick:
Noticing that many monomials share the same sequences of
$A_{i,q^k}^m$'s, we store the data into a tree so that we do not need
to repeat the common part.
By this trick, it becomes uncertain how much size we need in advance,
as we mentioned above.

Next let us turn to coloring.
By \cite{Na-qchar-main}, $\chi_{q,t}(L(P)) = \sum_m a_m(t) m$ is given
by the Poincar\'e polynomials of various graded quiver varieties
corresponding to $m$. Therefore the degree of the
coefficient $a_m(t)$ is equal to the (real) dimension of the variety
corresponding to $m$. On the other hand, the dimension of the graded
quiver variety is bounded by the half of the ordinary quiver variety
containing it. For the $5^{\mathrm{th}}$ fundamental representaion,
the maximum (among various connected components) of the dimension is
equal to $60$. Therefore the maximum of the degree is $30$. As
$a_{m,i}(t)$ is given by a virtual Hodge polynomial of a certain
stratum of the graded quiver variety, the degree is also less than or
equal to $30$. As $a_m(t)$, $a_{m,i}(t)$ are polynomials in $t^2$, we
have $30/2 + 1 = 16$ coefficients. Therefore we must record
$16\times 8$ integers for each monomial. We did not know how
large integers were in advance. As a result of our calculation, it
turns out we can store it into a \verb+short int+. Then we would need
$16\times 8 \times 16\mathrm{bit} = 256\mathrm{byte}$ for each
monomial. This is huge size, though it could be handled by our computer
probably.
However we note that many monomials $m$ have coefficient $a_m(t) = 1$. We
store $a_{m,i}(t)$ for those monomials in a special format to save the
size of data. As we do not need $a_{m,i}(t)$ for the final result,
they are not included.
(As a result of our calculation we find $4639565354$ among
$6899079264$ monomials have this property.)

We have explained the total size of the data so far. In practice, it
is more important to know how much memory is required in the course of
the calculation.
For the simplicity of the program, we replace the ordering $<$ among
monomials by more manageable ordering given by
\begin{equation*}
    \operatorname{depth}m \defeq \sum_{i,a} v_{i,a}(m,m_P).
\end{equation*}
Therefore the $l$-highest weight vector has depth $0$,
$Y_{5,1}A_{5,q}^{-1}$ has depth $1$, etc. We expand the monomial of
depth $0$, then monomials with depth $1$, monomials with depth $2$,
and so on. When we expand all monomials of given depth, we store all
obtained monomials together with coloring in memory. As a single
monomial appears many times in the expansions at various monomials, it
is not practical to save the data in the hard disk.
Therefore the most crucial point is to save the size of data so that
the program requires, in a fixed depth, up to $200\mathrm{Gbyte}$ of
memory, which is the limit of the supercomputer. We estimated
the memory requirement by that for $4^{\mathrm{th}}$ $l$--fundamental
representation as above, and we guessed that the calculation was
possible. This turns out to be true fortunately.

\section{Results}\label{sec:results}

We only consider the $5^{\mathrm{th}}$ $l$--fundamental representation
of $\Ul$ with $\g = E_8$.

As the final result is a huge polynomial, we cannot give it here. So
we only give a part of the information. The monomial whose
coefficient with the highest degree $t^{30}$ is
\begin{multline*}
   (1 + 4t^2 + 10t^4 + 20t^6 + 33t^8 + 47t^{10} + 59t^{12} + 66t^{14}
\\
   + 66t^{16} + 59t^{18} + 47t^{20} + 33t^{22} + 20t^{24} + 10t^{26} +
   4t^{28} + t^{30})
\\
   \times Y_{1,q^{14}} Y_{1,q^{16}}^{-1} Y_{3,q^{14}}^2
   Y_{3,q^{16}}^{-2} Y_{5,q^{14}}^3 Y_{5,q^{16}}^{-3} Y_{7,q^{14}} Y_{7,q^{16}}^{-1}.
\end{multline*}
The coefficient is the Poincar\'e polynomial of a certain graded
quiver variety.

We define the $t$--graded character by
\begin{equation*}
   \operatorname{ch}_t(L(P))
   = \left.\qch{q,t}(L(P))\right|_{Y_{i,a}\to y_i}.
\end{equation*}
If we put $t=1$, it becomes the ordinary character of the restriction
of $L(P)$ to $\Uq$. It is also equal to the generating function of the
Poincar\'e polynomials of the quiver varieties, where the degree $0$
corresponding to the middle degree. For example, the coefficient of the 
weight $0$ is
\begin{multline*}
1357104 + 2232771{}t^2 + 2002423{}t^4 + 1317308{}t^6 + 716312{}t^8 +
342421{}t^{10} + 148512{}t^{12}
\\
+ 59490{}t^{14} + 22162{}t^{16} + 7687{}t^{18} +
2463{}t^{20} + 726{}t^{22} + 192{}t^{24} + 44{}t^{26} + 8{}t^{28} + t^{30}.
\end{multline*}

Let $V(\lambda)$ denote the irreducible highest weight representation
of $\Uq$ with the highest weight $\lambda$. Let $\operatorname{ch}
V(\lambda)$ be its character. If we write
\begin{equation*}
   \operatorname{ch}_t L(P)
   = \sum_\lambda M(P,\lambda,t) \operatorname{ch} V(\lambda),
\end{equation*}
the coefficient $M(P,\lambda,t)$ is specialized to the multiplicity of
$V(\lambda)$ in the restriction of $L(P)$ at $t=1$. The fermionic
formula mentioned in the Introduction is a conjectural expression of
$M(P,\lambda,t)$ (for $P$ corresponding to the Kirillov-Reshetikhin
modules).

As we have computed $\qch{q,t}(L(P))$, $M(P,\lambda,t)$ can be given
if we compute $V(\lambda)$. Let us compute $V(\lambda)$ by the method in
\cite[7.1.1]{Na-qchar}, i.e.\ 
\begin{equation*}
   V(\lambda) = \left.\qch{q,t}(L(Q))\right|_{Y_{i,a}\to y_i, t\to 0},
\end{equation*}
where $Q$ corresponding to $\lambda$ is given as follows: We choose an
orientation for each edge of the Dynkin diagram and choose a function 
$m\colon I\to \Z$ such that $m(i) - m(j) = 1$ for an oriented
edge $i\to j$. Then we take 
\begin{equation*}
   Q_i(u) = (1 - uq^{m(i)})^{\langle\lambda, h_i\rangle}.
\end{equation*}
For this choice of $Q$, it is known that 
\(
   \operatorname{ch}_t(L(Q)) =
   \qch{q,t}(L(Q))|_{Y_{i,a}\to y_i}
\)
is equal to the generating function of shifted Poincar\'e polynomial
of the quiver variety as above. In particular, it is independent of
the choice of the orientation. For each dominant weight $\lambda$
appearing in $\operatorname{ch}_t L(P)$, we choose $Q = Q_\lambda$ as
above and define matrices
\(
   P(t) = (P_{\lambda\mu}(t))
\)
and
\(
   IC(t) = (IC_{\lambda\mu}(t))
\)
by
\begin{gather*}
   \operatorname{ch}_t L(Q_\lambda) 
   = \sum_\mu P_{\lambda\mu}(t) e^{\mu} + \text{non dominant terms},
\\
   \operatorname{ch}_t L(Q_\lambda) 
   = \sum_\mu IC_{\lambda\mu}(t) \operatorname{ch} V(\mu).
\end{gather*}
Then we have
\begin{equation*}
   IC(t) = P(t) P(0)^{-1}
\end{equation*}
By \cite{Na-qchar,Na-qchar-main} $IC_{\lambda\mu}(t)$ is the
Poincar\'e polynomial of the stalk of the intersection cohomology
sheaf of a stratum of the quiver variety corresponding to $\lambda$ at
a point in the stratum corresponding to $\mu$. In our case it is given
by
\begin{equation*}
   IC(t) =
\left(
\begin{tabular}{c|c|c}
   Table~\ref{tab:table1} & &
\\\cline{1-1}
   & \raisebox{1.5ex}[0cm][0cm]{Table~\ref{tab:table2}} & 
   \raisebox{1.5ex}[0cm][0cm]{Table~\ref{tab:table3}}
\\\cline{2-3}
   \raisebox{1.5ex}[0cm][0cm]{0} & 0 & Table~\ref{tab:table4}
\end{tabular}
\right),
\end{equation*}
where $y_i = e^{\varpi_i}$.
The first row gives $\operatorname{ch}_t(L(P))$ for the
$5^{\mathrm{th}}$ $l$--fundamental representation $L(P)$. We see that
it coincides with the conjectural formula in \cite{HKOTY}. The same
assertion for other $l$--fundamental representations can be proved by
invoking other rows. The same can be proved for types $E_6$, $E_7$ in
the same manner.

\begin{table}
\begin{sideways}
\newcolumntype{L}{>{$\scriptstyle}l<{$}}
\begin{tabular}{>{$}c<{$}|LLLLLLLLLL}
&
\multicolumn{1}{c}{$\vi$} & 
\multicolumn{1}{c}{$\vii$} &
\multicolumn{1}{c}{$\viii$} & 
\multicolumn{1}{c}{$\viv$} & 
\multicolumn{1}{c}{$\vv$} & 
\multicolumn{1}{c}{$\vvi$} & 
\multicolumn{1}{c}{$\vvii$} & 
\multicolumn{1}{c}{$\vviii$} &
\multicolumn{1}{c}{$\vix$} &
\multicolumn{1}{c}{$\vx$}
\\
\hline
\vi &
1&{t}^{2}&{t}^{2}+{t}^{4}&{t}^{4}&{
t}^{2}+{t}^{4}+{t}^{6}&{t}^{6}&{t}^{2}+2\,{t}^{4}+2\,{t}^{6}+{t}^{8}&2
\,{t}^{4}+{t}^{6}+{t}^{8}&{t}^{2}+2\,{t}^{4}+2\,{t}^{6}+{t}^{8}+{t}^{
10}&2\,{t}^{6}+{t}^{8}+{t}^{10}\\
\vii &
0&1&{t}^{2}&{t}^{2}
&{t}^{2}+{t}^{4}&{t}^{4}&{t}^{2}+2\,{t}^{4}+{t}^{6}&{t}^{2}+{t}^{4}+{t
}^{6}&{t}^{2}+2\,{t}^{4}+{t}^{6}+{t}^{8}&2\,{t}^{4}+{t}^{6}+{t}^{8}
\\
\viii &
0&0&1&0&{t}^{2}&{t}^{2}&{t}^{2}+{t}^{4}&{t}^{2}+{t
}^{4}&{t}^{2}+{t}^{4}+{t}^{6}&{t}^{4}+{t}^{6}\\
\viv &
0&0&0
&1&{t}^{2}&0&{t}^{2}+{t}^{4}&{t}^{4}&{t}^{4}+{t}^{6}&{t}^{2}+{t}^{4}+{
t}^{6}\\
\vv &
0&0&0&0&1&0&{t}^{2}&{t}^{2}&{t}^{2}+{t}^{4}&
{t}^{2}+{t}^{4}\\
\vvi &
0&0&0&0&0&1&0&{t}^{2}&{t}^{4}&{t}^{
4}\\
\vvii &
0&0&0&0&0&0&1&0&{t}^{2}&0\\
\vviii &
0&0
&0&0&0&0&0&1&{t}^{2}&{t}^{2}\\
\vix &
0&0&0&0&0&0&0&0&1&0
\\
\vx &
0&0&0&0&0&0&0&0&0&1
\end{tabular}
\end{sideways}
\caption{}
\label{tab:table1}
\end{table}

\begin{table}
\begin{sideways}
\newcolumntype{L}{>{$\scriptstyle}p{2.3cm}<{$}}
\begin{tabular}{>{$}c<{$}|LLLLLLLL}
&
\multicolumn{1}{c}{$\vxi$} & 
\multicolumn{1}{c}{$\vxii$} & 
\multicolumn{1}{c}{$\vxiii$} & 
\multicolumn{1}{c}{$\vxiv$} & 
\multicolumn{1}{c}{$\vxv$} & 
\multicolumn{1}{c}{$\vxvi$} &
\multicolumn{1}{c}{$\vxvii$} &
\multicolumn{1}{c}{$\vxviii$}
\\
\hline
%
\vi &
3\,{t}^{4}+4\,{t}^{6}+4\,{t}^{8}+2\,{
t}^{10}+{t}^{12}&2\,{t}^{4}+5\,{t}^{6}+5\,{t}^{8}+3\,{t}^{10}+2\,{t}^{
12}+{t}^{14}&3\,{t}^{6}+2\,{t}^{8}+3\,{t}^{10}+{t}^{12}+{t}^{14}&2\,{t
}^{4}+4\,{t}^{6}+6\,{t}^{8}+4\,{t}^{10}+3\,{t}^{12}+{t}^{14}+{t}^{16}&
{t}^{8}+{t}^{12}&2\,{t}^{6}+5\,{t}^{8}+5\,{t}^{10}+3\,{t}^{12}+2\,{t}^
{14}+{t}^{16}&5\,{t}^{6}+5\,{t}^{8}+7\,{t}^{10}+4\,{t}^{12}+3\,{t}^{14
}+{t}^{16}+{t}^{18}& 2\,{t}^{6}+9\,{t}^{8}+10\,{t}^{10}+10\,{t}^{12}
+6\,{t}^{14}+4\,{t}^{16}+2\,{t}^{18}+{t}^{20}
\\
\vii &
{t}^{2}
+3\,{t}^{4}+4\,{t}^{6}+2\,{t}^{8}+{t}^{10}&3\,{t}^{4}+5\,{t}^{6}+3\,{t
}^{8}+2\,{t}^{10}+{t}^{12}&2\,{t}^{4}+2\,{t}^{6}+3\,{t}^{8}+{t}^{10}+{
t}^{12}&2\,{t}^{4}+5\,{t}^{6}+4\,{t}^{8}+3\,{t}^{10}+{t}^{12}+{t}^{14}
&{t}^{6}+{t}^{10}&{t}^{4}+4\,{t}^{6}+5\,{t}^{8}+3\,{t}^{10}+2\,{t}^{12
}+{t}^{14}&2\,{t}^{4}+4\,{t}^{6}+7\,{t}^{8}+4\,{t}^{10}+3\,{t}^{12}+{t
}^{14}+{t}^{16}&6\,{t}^{6}+9\,{t}^{8}+10\,{t}^{10}+6\,{t}^{12}+4\,{t}^
{14}+2\,{t}^{16}+{t}^{18}\\
\viii &
{t}^{2}+3\,{t}^{4}+2\,{t}
^{6}+{t}^{8}&{t}^{2}+3\,{t}^{4}+3\,{t}^{6}+2\,{t}^{8}+{t}^{10}&{t}^{4}
+2\,{t}^{6}+{t}^{8}+{t}^{10}&3\,{t}^{4}+3\,{t}^{6}+3\,{t}^{8}+{t}^{10}
+{t}^{12}&{t}^{8}&2\,{t}^{4}+4\,{t}^{6}+3\,{t}^{8}+2\,{t}^{10}+{t}^{12
}&2\,{t}^{4}+5\,{t}^{6}+4\,{t}^{8}+3\,{t}^{10}+{t}^{12}+{t}^{14}&{t}^{
4}+6\,{t}^{6}+8\,{t}^{8}+6\,{t}^{10}+4\,{t}^{12}+2\,{t}^{14}+{t}^{16}
\\
\viv &
{t}^{2}+2\,{t}^{4}+2\,{t}^{6}+{t}^{8}&2\,{t}^{4}+3
\,{t}^{6}+2\,{t}^{8}+{t}^{10}&{t}^{2}+{t}^{4}+2\,{t}^{6}+{t}^{8}+{t}^{
10}&2\,{t}^{4}+3\,{t}^{6}+3\,{t}^{8}+{t}^{10}+{t}^{12}&{t}^{4}+{t}^{8}
&{t}^{4}+3\,{t}^{6}+3\,{t}^{8}+2\,{t}^{10}+{t}^{12}&4\,{t}^{6}+4\,{t}^
{8}+3\,{t}^{10}+{t}^{12}+{t}^{14}&2\,{t}^{4}+4\,{t}^{6}+8\,{t}^{8}+6\,
{t}^{10}+4\,{t}^{12}+2\,{t}^{14}+{t}^{16}\\
\vv &
{t}^{2}+2
\,{t}^{4}+{t}^{6}&{t}^{2}+3\,{t}^{4}+2\,{t}^{6}+{t}^{8}&{t}^{4}+{t}^{6
}+{t}^{8}&{t}^{2}+2\,{t}^{4}+3\,{t}^{6}+{t}^{8}+{t}^{10}&{t}^{6}&2\,{t
}^{4}+3\,{t}^{6}+2\,{t}^{8}+{t}^{10}&3\,{t}^{4}+4\,{t}^{6}+3\,{t}^{8}+
{t}^{10}+{t}^{12}&2\,{t}^{4}+6\,{t}^{6}+6\,{t}^{8}+4\,{t}^{10}+2\,{t}^
{12}+{t}^{14}\\
\vvi &
{t}^{2}+{t}^{4}+{t}^{6}&{t}^{4}+2\,{t
}^{6}+{t}^{8}&{t}^{4}+{t}^{8}&2\,{t}^{6}+{t}^{8}+{t}^{10}&{t}^{6}&{t}^
{2}+2\,{t}^{4}+2\,{t}^{6}+2\,{t}^{8}+{t}^{10}&2\,{t}^{4}+2\,{t}^{6}+3
\,{t}^{8}+{t}^{10}+{t}^{12}&{t}^{4}+4\,{t}^{6}+4\,{t}^{8}+4\,{t}^{10}+
2\,{t}^{12}+{t}^{14}\\
\vvii &
{t}^{2}+{t}^{4}&{t}^{2}+{t}^{4
}+{t}^{6}&{t}^{2}+{t}^{4}+{t}^{6}&{t}^{2}+2\,{t}^{4}+{t}^{6}+{t}^{8}&0
&{t}^{4}+{t}^{6}+{t}^{8}&2\,{t}^{4}+2\,{t}^{6}+{t}^{8}+{t}^{10}&2\,{t}
^{4}+4\,{t}^{6}+3\,{t}^{8}+2\,{t}^{10}+{t}^{12}
\\
\vviii &
{t}^{2}+{t}^{4}&{t}^{2}+2\,{t}^{4}+{t}^{6}&{t}^{6}&2\,{t}^{4}+{t}^{6}+{t}
^{8}&{t}^{4}&{t}^{2}+2\,{t}^{4}+2\,{t}^{6}+{t}^{8}&{t}^{2}+2\,{t}^{4}+
3\,{t}^{6}+{t}^{8}+{t}^{10}&3\,{t}^{4}+4\,{t}^{6}+4\,{t}^{8}+2\,{t}^{
10}+{t}^{12}
\\
\vix &
{t}^{2}&{t}^{2}+{t}^{4}&{t}^{4}&{t}^{2
}+{t}^{4}+{t}^{6}&0&{t}^{4}+{t}^{6}&{t}^{2}+2\,{t}^{4}+{t}^{6}+{t}^{8}
&2\,{t}^{4}+3\,{t}^{6}+2\,{t}^{8}+{t}^{10}\\
\vx &
{t}^{2}&
{t}^{2}+{t}^{4}&{t}^{4}&{t}^{4}+{t}^{6}&{t}^{2}&{t}^{2}+2\,{t}^{4}+{t}
^{6}&2\,{t}^{4}+{t}^{6}+{t}^{8}&{t}^{2}+2\,{t}^{4}+4\,{t}^{6}+2\,{t}^{
8}+{t}^{10}
\\
\vxi &
1&{t}^{2}&{t}^{2}&{t}^{2}+{t}^{4}&0&{t}
^{2}+{t}^{4}&{t}^{2}+{t}^{4}+{t}^{6}&{t}^{2}+3\,{t}^{4}+2\,{t}^{6}+{t}
^{8}\\
\vxii &
0&1&0&{t}^{2}&0&{t}^{2}&{t}^{2}+{t}^{4}&{t}^{2
}+2\,{t}^{4}+{t}^{6}\\
\vxiii &
0&0&1&{t}^{2}&0&0&{t}^{4}&{t}^
{2}+{t}^{4}+{t}^{6}\\
\vxiv &
0&0&0&1&0&0&{t}^{2}&{t}^{2}+{t}
^{4}\\
\vxv &
0&0&0&0&1&{t}^{2}+{t}^{4}&{t}^{6}&{t}^{4}+{t}^
{6}+{t}^{8}\\
\vxvi &
0&0&0&0&0&1&{t}^{2}&{t}^{2}+{t}^{4}
\\
\vxvii &
0&0&0&0&0&0&1&{t}^{2}\\
\vxviii&
0&0&0&0&0&0&0&1
\end{tabular}
\end{sideways}
\caption{}
\label{tab:table2}
\end{table}

\begin{table}
\begin{sideways}
\newcolumntype{L}{>{$\scriptstyle}p{3.3cm}<{$}}
\begin{tabular}{>{$}c<{$}|LLLLLL}
&
\multicolumn{1}{c}{$\vxix$} & 
\multicolumn{1}{c}{$\vxx$} & 
\multicolumn{1}{c}{$\vxxi$} & 
\multicolumn{1}{c}{$\vxxii$} & 
\multicolumn{1}{c}{$\vxxiii$} & 
\multicolumn{1}{c}{$\vxxiv$}
\\
\hline
%
\vi &
{t}^{6}+5\,{t}^{8}+8\,{t}^{10}+7\,{t}^{
12}+6\,{t}^{14}+4\,{t}^{16}+2\,{t}^{18}+{t}^{20}+{t}^{22}&5\,{t}^{10}+
4\,{t}^{12}+6\,{t}^{14}+3\,{t}^{16}+3\,{t}^{18}+{t}^{20}+{t}^{22}&3\,{
t}^{8}+6\,{t}^{10}+11\,{t}^{12}+8\,{t}^{14}+7\,{t}^{16}+4\,{t}^{18}+3
\,{t}^{20}+{t}^{22}+{t}^{24}&5\,{t}^{10}+6\,{t}^{12}+9\,{t}^{14}+6\,{t
}^{16}+6\,{t}^{18}+3\,{t}^{20}+2\,{t}^{22}+{t}^{24}+{t}^{26}&4\,{t}^{
12}+5\,{t}^{14}+8\,{t}^{16}+5\,{t}^{18}+5\,{t}^{20}+3\,{t}^{22}+2\,{t}
^{24}+{t}^{26}+{t}^{28}&{t}^{14}+3\,{t}^{18}+{t}^{20}+2\,{t}^{22}+{t}^
{24}+{t}^{26}+{t}^{30}
\\
\vii &
2\,{t}^{6}+7\,{t}^{8}+7\,{t}
^{10}+6\,{t}^{12}+4\,{t}^{14}+2\,{t}^{16}+{t}^{18}+{t}^{20}&4\,{t}^{8}
+4\,{t}^{10}+6\,{t}^{12}+3\,{t}^{14}+3\,{t}^{16}+{t}^{18}+{t}^{20}&{t}
^{6}+4\,{t}^{8}+10\,{t}^{10}+8\,{t}^{12}+7\,{t}^{14}+4\,{t}^{16}+3\,{t
}^{18}+{t}^{20}+{t}^{22}&3\,{t}^{8}+5\,{t}^{10}+9\,{t}^{12}+6\,{t}^{14
}+6\,{t}^{16}+3\,{t}^{18}+2\,{t}^{20}+{t}^{22}+{t}^{24}&3\,{t}^{10}+4
\,{t}^{12}+8\,{t}^{14}+5\,{t}^{16}+5\,{t}^{18}+3\,{t}^{20}+2\,{t}^{22}
+{t}^{24}+{t}^{26}&{t}^{12}+3\,{t}^{16}+{t}^{18}+2\,{t}^{20}+{t}^{22}+
{t}^{24}+{t}^{28}
\\
\viii &
4\,{t}^{6}+5\,{t}^{8}+6\,{t}^{10}
+4\,{t}^{12}+2\,{t}^{14}+{t}^{16}+{t}^{18}&{t}^{6}+3\,{t}^{8}+5\,{t}^{
10}+3\,{t}^{12}+3\,{t}^{14}+{t}^{16}+{t}^{18}&{t}^{6}+7\,{t}^{8}+7\,{t
}^{10}+7\,{t}^{12}+4\,{t}^{14}+3\,{t}^{16}+{t}^{18}+{t}^{20}&3\,{t}^{8
}+6\,{t}^{10}+6\,{t}^{12}+6\,{t}^{14}+3\,{t}^{16}+2\,{t}^{18}+{t}^{20}
+{t}^{22}&3\,{t}^{10}+6\,{t}^{12}+5\,{t}^{14}+5\,{t}^{16}+3\,{t}^{18}+
2\,{t}^{20}+{t}^{22}+{t}^{24}&2\,{t}^{14}+{t}^{16}+2\,{t}^{18}+{t}^{20
}+{t}^{22}+{t}^{26}
\\
\viv &
2\,{t}^{6}+5\,{t}^{8}+6\,{t}^{
10}+4\,{t}^{12}+2\,{t}^{14}+{t}^{16}+{t}^{18}&2\,{t}^{6}+2\,{t}^{8}+5
\,{t}^{10}+3\,{t}^{12}+3\,{t}^{14}+{t}^{16}+{t}^{18}&{t}^{6}+5\,{t}^{8
}+7\,{t}^{10}+7\,{t}^{12}+4\,{t}^{14}+3\,{t}^{16}+{t}^{18}+{t}^{20}&{t
}^{6}+2\,{t}^{8}+6\,{t}^{10}+6\,{t}^{12}+6\,{t}^{14}+3\,{t}^{16}+2\,{t
}^{18}+{t}^{20}+{t}^{22}&2\,{t}^{8}+2\,{t}^{10}+6\,{t}^{12}+5\,{t}^{14
}+5\,{t}^{16}+3\,{t}^{18}+2\,{t}^{20}+{t}^{22}+{t}^{24}&{t}^{10}+2\,{t
}^{14}+{t}^{16}+2\,{t}^{18}+{t}^{20}+{t}^{22}+{t}^{26}
\\
\vv &
{t}^{4}+4\,{t}^{6}+6\,{t}^{8}+4\,{t}^{10}+2\,{t}^{
12}+{t}^{14}+{t}^{16}&{t}^{6}+4\,{t}^{8}+3\,{t}^{10}+3\,{t}^{12}+{t}^{
14}+{t}^{16}&3\,{t}^{6}+6\,{t}^{8}+7\,{t}^{10}+4\,{t}^{12}+3\,{t}^{14}
+{t}^{16}+{t}^{18}&{t}^{6}+4\,{t}^{8}+6\,{t}^{10}+6\,{t}^{12}+3\,{t}^{
14}+2\,{t}^{16}+{t}^{18}+{t}^{20}&{t}^{8}+4\,{t}^{10}+5\,{t}^{12}+5\,{
t}^{14}+3\,{t}^{16}+2\,{t}^{18}+{t}^{20}+{t}^{22}&{t}^{12}+{t}^{14}+2
\,{t}^{16}+{t}^{18}+{t}^{20}+{t}^{24}
\\
\vvi &
{t}^{6}+4\,{t
}^{8}+4\,{t}^{10}+2\,{t}^{12}+{t}^{14}+{t}^{16}&{t}^{4}+{t}^{6}+4\,{t}
^{8}+2\,{t}^{10}+3\,{t}^{12}+{t}^{14}+{t}^{16}&3\,{t}^{6}+3\,{t}^{8}+6
\,{t}^{10}+4\,{t}^{12}+3\,{t}^{14}+{t}^{16}+{t}^{18}&3\,{t}^{8}+3\,{t}
^{10}+6\,{t}^{12}+3\,{t}^{14}+2\,{t}^{16}+{t}^{18}+{t}^{20}&{t}^{8}+4
\,{t}^{10}+3\,{t}^{12}+5\,{t}^{14}+3\,{t}^{16}+2\,{t}^{18}+{t}^{20}+{t
}^{22}&2\,{t}^{12}+2\,{t}^{16}+{t}^{18}+{t}^{20}+{t}^{24}
\\
\vvii &
{t}^{4}+3\,{t}^{6}+3\,{t}^{8}+2\,{t}^{10}+{t}^{12}
+{t}^{14}&{t}^{6}+2\,{t}^{8}+2\,{t}^{10}+{t}^{12}+{t}^{14}&2\,{t}^{6}+
4\,{t}^{8}+3\,{t}^{10}+3\,{t}^{12}+{t}^{14}+{t}^{16}&{t}^{6}+3\,{t}^{8
}+4\,{t}^{10}+3\,{t}^{12}+2\,{t}^{14}+{t}^{16}+{t}^{18}&{t}^{8}+3\,{t}
^{10}+3\,{t}^{12}+3\,{t}^{14}+2\,{t}^{16}+{t}^{18}+{t}^{20}&{t}^{12}+{
t}^{14}+{t}^{16}+{t}^{18}+{t}^{22}
\\
\vviii &
{t}^{4}+4\,{t}^{
6}+4\,{t}^{8}+2\,{t}^{10}+{t}^{12}+{t}^{14}&3\,{t}^{6}+2\,{t}^{8}+3\,{
t}^{10}+{t}^{12}+{t}^{14}&{t}^{4}+3\,{t}^{6}+6\,{t}^{8}+4\,{t}^{10}+3
\,{t}^{12}+{t}^{14}+{t}^{16}&2\,{t}^{6}+3\,{t}^{8}+6\,{t}^{10}+3\,{t}^
{12}+2\,{t}^{14}+{t}^{16}+{t}^{18}&3\,{t}^{8}+3\,{t}^{10}+5\,{t}^{12}+
3\,{t}^{14}+2\,{t}^{16}+{t}^{18}+{t}^{20}&{t}^{10}+2\,{t}^{14}+{t}^{16
}+{t}^{18}+{t}^{22}
\\
\vix &
2\,{t}^{4}+3\,{t}^{6}+2\,{t}^{8
}+{t}^{10}+{t}^{12}&{t}^{6}+2\,{t}^{8}+{t}^{10}+{t}^{12}&3\,{t}^{6}+3
\,{t}^{8}+3\,{t}^{10}+{t}^{12}+{t}^{14}&{t}^{6}+4\,{t}^{8}+3\,{t}^{10}
+2\,{t}^{12}+{t}^{14}+{t}^{16}&{t}^{8}+3\,{t}^{10}+3\,{t}^{12}+2\,{t}^
{14}+{t}^{16}+{t}^{18}&{t}^{12}+{t}^{14}+{t}^{16}+{t}^{20}
\\
\vx &
{t}^{4}+3\,{t}^{6}+2\,{t}^{8}+{t}^{10}+{t}^{12}&2
\,{t}^{4}+2\,{t}^{6}+3\,{t}^{8}+{t}^{10}+{t}^{12}&{t}^{4}+4\,{t}^{6}+4
\,{t}^{8}+3\,{t}^{10}+{t}^{12}+{t}^{14}&{t}^{4}+{t}^{6}+5\,{t}^{8}+3\,
{t}^{10}+2\,{t}^{12}+{t}^{14}+{t}^{16}&2\,{t}^{6}+2\,{t}^{8}+5\,{t}^{
10}+3\,{t}^{12}+2\,{t}^{14}+{t}^{16}+{t}^{18}&{t}^{8}+2\,{t}^{12}+{t}^
{14}+{t}^{16}+{t}^{20}
\\
\vxi &
2\,{t}^{4}+2\,{t}^{6}+{t}^{8
}+{t}^{10}&{t}^{4}+2\,{t}^{6}+{t}^{8}+{t}^{10}&2\,{t}^{4}+3\,{t}^{6}+3
\,{t}^{8}+{t}^{10}+{t}^{12}&3\,{t}^{6}+3\,{t}^{8}+2\,{t}^{10}+{t}^{12}
+{t}^{14}&{t}^{6}+3\,{t}^{8}+3\,{t}^{10}+2\,{t}^{12}+{t}^{14}+{t}^{16}
&{t}^{10}+{t}^{12}+{t}^{14}+{t}^{18}
\\
\vxii &
{t}^{2}+2\,{t}
^{4}+{t}^{6}+{t}^{8}&{t}^{4}+{t}^{6}+{t}^{8}&2\,{t}^{4}+3\,{t}^{6}+{t}
^{8}+{t}^{10}&{t}^{4}+3\,{t}^{6}+2\,{t}^{8}+{t}^{10}+{t}^{12}&{t}^{6}+
3\,{t}^{8}+2\,{t}^{10}+{t}^{12}+{t}^{14}&{t}^{10}+{t}^{12}+{t}^{16}
\\
\vxiii &
{t}^{4}+{t}^{6}+{t}^{8}&{t}^{4}+{t}^{8}&2\,{t}^{6}
+{t}^{8}+{t}^{10}&{t}^{4}+{t}^{6}+2\,{t}^{8}+{t}^{10}+{t}^{12}&{t}^{6}
+{t}^{8}+2\,{t}^{10}+{t}^{12}+{t}^{14}&{t}^{8}+{t}^{12}+{t}^{16}
\\
\vxiv &
{t}^{2}+{t}^{4}+{t}^{6}&{t}^{6}&2\,{t}^{4}+{t}^{6}
+{t}^{8}&{t}^{4}+2\,{t}^{6}+{t}^{8}+{t}^{10}&{t}^{6}+2\,{t}^{8}+{t}^{
10}+{t}^{12}&{t}^{10}+{t}^{14}
\\
\vxv &
{t}^{8}+{t}^{10}&{t}^{2}+{t}^{4}+2\,{t}
^{6}+{t}^{8}+{t}^{10}&{t}^{4}+2\,{t}^{6}+2\,{t}^{8}+{t}^{10}+{t}^{12}&
{t}^{6}+{t}^{8}+2\,{t}^{10}+{t}^{12}+{t}^{14}&{t}^{4}+{t}^{6}+3\,{t}^{
8}+2\,{t}^{10}+2\,{t}^{12}+{t}^{14}+{t}^{16}&{t}^{6}+{t}^{10}+{t}^{12}
+{t}^{14}+{t}^{18}
\\
\vxvi &
{t}^{4}+{t}^{6}&{t}^{2}+{t}^{4}+
{t}^{6}&{t}^{2}+2\,{t}^{4}+{t}^{6}+{t}^{8}&{t}^{4}+2\,{t}^{6}+{t}^{8}+
{t}^{10}&{t}^{4}+2\,{t}^{6}+2\,{t}^{8}+{t}^{10}+{t}^{12}&{t}^{8}+{t}^{
10}+{t}^{14}
\\
\vxvii &
{t}^{2}+{t}^{4}&{t}^{4}&{t}^{2}+{t}^{4
}+{t}^{6}&2\,{t}^{4}+{t}^{6}+{t}^{8}&2\,{t}^{6}+{t}^{8}+{t}^{10}&{t}^{
8}+{t}^{12}
\\
\vxviii&
{t}^{2}&{t}^{2}&{t}^{2}+{t}^{4}&{t}^{2}
+{t}^{4}+{t}^{6}&2\,{t}^{4}+{t}^{6}+{t}^{8}&{t}^{6}+{t}^{10}
\end{tabular}
\end{sideways}
\caption{}
\label{tab:table3}
\end{table}

\begin{table}
\newcolumntype{L}{>{$\scriptstyle}l<{$}}
\begin{tabular}{>{$}c<{$}|LLLLLL}
&
\multicolumn{1}{c}{$\vxix$} & 
\multicolumn{1}{c}{$\vxx$} & 
\multicolumn{1}{c}{$\vxxi$} & 
\multicolumn{1}{c}{$\vxxii$} & 
\multicolumn{1}{c}{$\vxxiii$} & 
\multicolumn{1}{c}{$\vxxiv$}
\\
\hline
%
\vxix&
1&0&{t}^{2}&{t}^{2}+{t}^{4}&{t}^{4}+{t}^{6}&{t}^{8
}
\\
\vxx&
0&1&{t}^{2}&{t}^{4}&{t}^{2}+{t}^{4}+{t}^{6}&{t}^{
4}+{t}^{8}
\\
\vxxi&
0&0&1&{t}^{2}&{t}^{2}+{t}^{4}&{t}^{6}
\\
\vxxii&
0&0&0&1&{t}^{2}&{t}^{4}
\\
\vxxiii&
0&0&0&0&1&{t}^{2}
\\
\vxxiv&
0&0&0&0&0&1
\end{tabular}
\caption{}
\label{tab:table4}
\end{table}

\end{document}